\input ./style/arxiv-ba.cfg
\documentclass[ba,linksfromyear,preprint]{imsart}
\makeatletter
   \@ifpackageloaded{natbib}{}{\usepackage{natbib}}
\makeatother

\pubyear{2015}
\volume{10}
\issue{3}
\firstpage{523}
\lastpage{552}
\doi{10.1214/14-BA922}

\newcommand{\bs}{\boldsymbol}
\newcommand{\lik}{\mathbb L}

\newtheorem{theorem}{Theorem}
\newtheorem{coroll}{Corollary}

\begin{document}

\begin{frontmatter}
\title{Generalized Quantile Treatment Effect: A~Flexible Bayesian
Approach Using Quantile Ratio Smoothing}
\runtitle{Generalized Quantile Treatment Effect}

\begin{aug}
\author[a]{\fnms{Sergio} \snm{Venturini}\ead
[label=e1]{sergio.venturini@unibocconi.it}},
\author[b]{\fnms{Francesca} \snm{Dominici}\ead
[label=e2]{fdominic@hsph.harvard.edu}},
\and
\author[c]{\fnms{Giovanni} \snm{Parmigiani}\ead
[label=e3]{gp@jimmy.harvard.edu}}

\runauthor{S. Venturini, F. Dominici and G. Parmigiani}


\address[a]{CERGAS, Universit\`a Bocconi, Via R\"ontgen 1, 20136 Milano,
Italy, \printead{e1}}
\address[b]{Department of Biostatistics, Harvard School of Public Health,
655 Huntington Ave, Boston MA 02115, USA, \printead{e2}}
\address[c]{Department of Biostatistics, Harvard School of Public Health
and Department of Biostatistics and Computational Biology, Dana Farber Cancer
Institute, 44 Binney Street, Boston MA 02115, USA, \printead{e3}}
\end{aug}

%
\begin{abstract}
We propose a new general approach for estimating the effect of a binary
treatment on a continuous and potentially highly skewed response variable,
the \textit{generalized quantile treatment effect} (GQTE). The GQTE is defined
as the difference between a function of the quantiles under the two treatment
conditions. As such, it represents a generalization over the standard approaches
typically used for estimating a treatment effect (i.e., the average treatment
effect and the quantile treatment effect) because it allows the
comparison of
any arbitrary characteristic of the outcome's distribution under the
two treatments.
Following \citet{DomCopeZeger}, we assume that a pre-specified transformation
of the two quantiles is modeled as a smooth function of the
percentiles. This
assumption allows us to link the two quantile functions and thus to borrow
information from one distribution to the other. The main theoretical
contribution
we provide is the analytical derivation of a closed form expression for the
likelihood of the model. Exploiting this result we propose a novel Bayesian
inferential methodology for the GQTE.
We show some finite sample properties of our approach through a simulation
study which confirms that in some cases it performs better than other
nonparametric methods. As an illustration we finally apply our methodology
to the $1987$ National Medicare Expenditure Survey data to estimate the
difference in the single hospitalization medical cost distributions
between cases (i.e., subjects affected by smoking attributable
diseases) and controls.
\end{abstract}

%


\end{frontmatter}


\section{Introduction}\label{introduction}
The effect of a treatment on an outcome is often the main parameter of
interest in many scientific fields. 
The standard approach used to estimate it is the so called average
treatment effect (ATE), the difference between the expected values of
the response's distributions under the two treatment regimes. While
intuitive and useful in many situations, it suffers from some
limitations; in particular, it becomes highly biased when the response
is skewed. 

A further drawback of the ATE is its coarseness as a summary of the
distance between the expected value of the response's distributions
under the two treatments. It is a matter of fact indeed that the effect
of the treatment on the outcome often varies as we move from the lower
to the upper tail of the outcome's distribution. This limitation of the
ATE has been addressed in the literature by introducing the so called
quantile treatment effect (QTE), the difference between the response's
distribution quantiles under the two treatments
\citep{Abadieetal2002,ChernozhukovHansen,Firpo2007,FrolichMelly2008}. 

In this paper we propose a more general measure of the effect of a
binary treatment on a continuous outcome. We call it the \emph
{generalized quantile treatment effect} (GQTE), defined as
%
\begin{equation}\label{intro1}
\Delta_{g}(p) = g(Q_1(p)) - g(Q_2(p)),
\end{equation}
where $Q_{1}(p)$ and $Q_{2}(p)$ represent the quantile functions of the
outcome under the two treatment conditions and $g(\cdot)$ is an
arbitrary but known function of the quantiles. For example, if $g(\cdot
)$ is chosen to be the identity function, then the GQTE simplifies to
the QTE, while if $g(\cdot)$ is the integral over the percentile $p$,
the GQTE becomes equivalent to the ATE. The GQTE is a new parameter
which generalizes the existing approaches for estimating a treatment effect.

To estimate and formulate inferences about the GQTE we propose a
Bayesian approach that can accommodate both symmetric and skewed
outcomes, as well as situations where the sample size under a treatment
condition (cases) is much smaller than the sample size under the other
treatment condition (controls). In particular, we assume
%
\begin{equation}\label{intro2}
h\left(\frac{Q_{1}(p)}{Q_{2}(p)}\right) = s(p)\,,
\end{equation}
where $h$ is a monotone function and $s$ is assumed to be smooth. In
other words, we assume that the transformed quantile ratio is a smooth
function of the percentile $p$. The idea of smoothly modeling the ratio
of the quantiles has been first introduced by \citet{DomCopeZeger},
who exploited it by proposing a nonparametric estimator of the mean
difference between two populations. Here we generalize their approach
by permitting the comparison of any characteristic of the outcome's
distributions under the two treatments.


An important theoretical contribution of this paper is the derivation
of a closed form expression for the model likelihood. We show that it
is possible to obtain an analytically tractable form for the $Y_{2}$
density (the controls) without explicitly specifying a model for it.
Clearly the likelihood is needed to carry out the Bayesian estimation
but in principle it could be employed for classic likelihood procedures
as well. Moreover, our proposed approach allows one to borrow strength
from one sample to the other, thus improving efficiency in the
estimation of the quantiles \citep{DomCopeZeger}.


As an illustration, we apply our method to the comparison of the single
hospitalization medical costs distribution between subjects with and
without smoking attributable diseases. The data set we use is the
National Medical Expenditures Survey (NMES) supplemented by the Adult
Self-Administered Questionnaire Household Survey.

The paper is organized as follows. In Section \ref{GQTE} we define the
new parameter $\Delta_g(p)$ and illustrate some quantile-based
measures that will be used in the paper. In Section~\ref{estimation}
we provide details of the estimation approach together with some
special cases. We then present the results of a simulation study in
Section \ref{simstudy} through which we conclude that under a broad
set of conditions our approach performs better than other flexible
methods for comparing two distributions. In Section \ref
{data_analysis} we illustrate the results of the data analysis on the
NMES data set. Section \ref{discussion} concludes the paper with a
discussion and some final remarks.

\section{The Generalized Quantile Treatment Effect (GQTE)}\label{GQTE}
Consider two positive continuous random variables $Y_1$ and $Y_2$ with
quantile functions $Q_1$ and $Q_2$, where
\[
Q_\ell(p) \equiv F^{-1}_\ell(p) \equiv\inf\{y: F_\ell(y) \geq p\}
\]
for $0<p<1$ and $\ell=1,2$. To compare $F_1$ and $F_2$ as flexibly as
possible we introduce the generalized quantile treatment effect, which
is defined as
%
\begin{equation}\label{GQTE_def}
\Delta_{g}(p) = g(Q_1(p)) - g(Q_2(p)),
\end{equation}
where $g(\cdot)$ is a known function of the quantiles. Notice that no
a priori assumptions are made about the admissible functions $g(\cdot
)$, thus potentially any function of the quantiles can be used.
Therefore, the GQTE provides a general approach to compare the
response's distributions under the two treatments. More precisely, by
properly choosing the function $g(\cdot)$, we can recover any specific
characteristic of the outcome's distributions $F_{1}$ and $F_{2}$ and,
through \eqref{GQTE_def}, their difference.

The simplest case arises when $g(x)=x$. In this case the GQTE
simplifies to
%
\begin{equation}\label{gqte_QTE}
\Delta(p)=Q_1(p)-Q_2(p),
\end{equation}
the so called (unconditional) QTE \citep{FrolichMelly2008}, sometimes
also named the percentile-specific effect between two populations
\citep[see for example][]{dometal06A,dometal07B}.

A second example is obtained by choosing $g(x)=\int{x}\,dp$, which produces
%
\begin{equation}\label{gqte_ATE}
\Delta=\int_0^1{Q_1(p)}\,dp-\int_0^1{Q_2(p)}\,dp \,,
\end{equation}
the extensively used ATE \citep[see for example][Chapter 21]{Wooldridge2010}.

These examples illustrate how the GQTE reduces to the two most used
parameters of interest for estimating a treatment effect, the ATE and
QTE. However, the GQTE can provide a variety of other useful measures.
In Appendix 1 
we illustrate some other interesting cases that usually are not taken
into consideration in the literature.

\section{Estimation Methodology}\label{estimation}
In this section we illustrate the procedure we developed for estimating
the GQTE. Our proposed approach is sufficiently general that it can be
used for any choice of $g(\cdot)$.

\subsection{Definitions and Model Assumptions}\label{bsq_def}
We assume that $Y_1 | \bs{\eta} \sim F_1(\cdot\,;\bs{\eta})$,
where $F_{1}$ is a given probability distribution depending upon a
vector of unknown parameters $\bs{\eta}$. For example, in the
application presented in Section \ref{data_analysis} we choose $F_{1}$
as a mixture distribution. To borrow information from one distribution
to the other, we assume that the transformed quantile ratio is a smooth
function of the percentiles with $\lambda$ degrees of freedom, that is
%
\begin{equation}\label{gqte_eq}
h\left( \frac{Q_1(p)}{Q_2(p)} \right)=s(p,\lambda)\,,\qquad0<p<1\,.
\end{equation}
The function $h(\cdot)$ is assumed to be monotone differentiable. It
represents a kind of link function and it is used to transform the
quantile ratio to account for the potential skewness of the $F_{1}$ and
$F_{2}$. The typical choice for skewed data is $h(x)=\log x$, while for
symmetric data distributions the identity function is the most
reasonable option.

For the sake of simplicity, we henceforth indicate the smooth function
$s(p,\lambda)$ with reference to the corresponding design matrix
$X(p,\lambda)$, so that it can be written as $X(p,\lambda)\bs{\beta
}$, where $\bs{\beta}$ is a vector of unknown parameters. More
explicitly, we assume that $s(p,\lambda) \equiv X(p,\lambda)\bs
{\beta}=\sum^{\lambda}_{k=0}{X_k(p)\beta_k}$, where $X_k(p)$ are
orthonormal basis functions with $X_0(p) = 1$. The number of degrees of
freedom $\lambda$ is a further parameter that has either to be chosen
or estimated from the data. In Subsection \ref{dfselection} we propose
a simple approach for eliciting it. The basis functions are usually
either splines or polynomials.

The main justification for assuming \eqref{gqte_eq} is that it allows
one to borrow information from both the response's distributions under
the two treatment conditions when we estimate the GQTE $\Delta
_{g}(p)$. Assumption \eqref{gqte_eq}, in fact, implies
%
\begin{equation}\label{gqte1}
Q_1(p)=Q_2(p)\,h^{-1}\left[X(p,\lambda)\,\bs{\beta}\right]\,,
\end{equation}
and also
%
\begin{equation}\label{gqte2}
Q_2(p)=Q_1(p) \left\{h^{-1}\left[X(p,\lambda)\,\bs{\beta}\right
]\right\}^{-1}\,,
\end{equation}
which, once substituted in \eqref{GQTE_def}, return
\[
\Delta_{g}(p) = g\Big( Q_2(p)\,h^{-1}\left[X(p,\lambda)\,\bs{\beta
}\right] \Big) - g\Big( Q_1(p) \left\{h^{-1}\left[X(p,\lambda)\,
\bs{\beta}\right]\right\}^{-1} \Big).
\]
For the special case where $g(x)=x$ and $h(x)=\log(x)$, \cite{DomCopeZeger}
have shown that under assumption \eqref{gqte_eq} it is
possible to obtain a more efficient estimator of $\Delta$ than the
sample mean difference and the maximum likelihood estimator assuming
that $Y_{1}$ and $Y_{2}$ are both log-normal.

Notice that, since the main interest in the paper resides in the
estimation of $\Delta_{g}(p)$, for which only $\bs{\beta}$ is
required, $\bs{\eta}$ is treated as nuisance (see Subsection \ref
{gqte_inference} for further details).

In this paper we propose a Bayesian approach for estimating $\Delta
_{g}(p)$ for any choice of $g(\cdot)$ and $h(\cdot)$. An interesting
feature of our estimation procedure for $\Delta_{g}(p)$, is that we
only need to specify the distribution function for $Y_{1}$. The
specification of $F_1$ together with the relationship \eqref{gqte_eq}
automatically determines a distributional assumption for $Y_2$. We
refer to the distribution of $Y_{2}$ induced by $F_{1}$ and assumption
\eqref{gqte_eq} as $F_2(\cdot\,;\bs{\beta}, \bs{\eta})$.

As a last remark for this section, we want to highlight the difference
between the function $g(\cdot)$, introduced in the previous section,
and $h(\cdot)$, defined above in \eqref{gqte_eq}. They should not be
confused because they have distinct roles: the former identifies the
response's characteristic we want to estimate for assessing the
treatment effect, while the latter has been introduced as a mechanism
to attenuate the possible skewness present in the data.

\subsection{Estimation Approach and Likelihood}\label{likapprox}
The steps involved in our estimation approach are summarized as follows:
\begin{enumerate}
\item Choose a (possibly flexible) density $f_{1}(y_{1} | \bs{\eta})$
for $Y_{1}$, a smoothing function $s(p,\lambda)$ (usually a spline or
a polynomial) and a value for $\lambda$;
\item From \eqref{gqte2} derive the density function of $Y_{2}$, that
we denote as $f_{2}(y_{2}|\bs{\beta},\bs{\eta})$. Note that, as
proved by Theorem \ref{theorem1} below, this density will depend on
the model parameter $\bs{\beta}$ as well as on the parameter $\bs
{\eta}$ through the $Y_{1}$ density.
\item Calculate the joint likelihood $\lik\left(\bs{\beta},\bs
{\eta}|\bs{y}_{1},\bs{y}_{2}\right)$ to use for finding the
posterior distribution of $\left( \bs{\beta}, \bs{\eta} \right)$
in a Markov Chain Monte Carlo (MCMC) algorithm.
\item Obtain the posterior distribution of any special case of the GQTE.
\end{enumerate}
The critical step in this sequence is represented by the calculation of
the likelihood, which we now describe.

Consider two i.i.d. samples $\left( y_{11},\ldots,y_{1n_1} \right)$
and $\left( y_{21},\ldots,y_{2n_2} \right)$ drawn independently from
the two populations $F_1(\cdot\,;\bs{\eta})$ and $F_2(\cdot\,;\bs
{\beta}, \bs{\eta})$. We refer to the former as the cases (or the
treated) and to the latter as the controls (or the untreated). We
assume that these distribution functions have densities $f_{1}(\cdot\,
;\bs{\eta})$ and $f_{2}(\cdot\,;\bs{\beta}, \bs{\eta})$
respectively. The likelihood function for our model is then given by
%
\begin{equation}\label{bsq_lik}
\lik\left(\bs{\beta},\bs{\eta}|y_{11},\ldots
,y_{1n_1},y_{21},\ldots,y_{2n_2}\right) = \prod
^{n_1}_{i=1}{f_1(y_{1i}|\bs{\eta})} \times\prod
^{n_2}_{j=1}{f_2(y_{2j}|\bs{\beta},\bs{\eta})} .
\end{equation}

Since we didn't state any specific distributional assumption for
$Y_{2}$, in principle we could not calculate the likelihood because we
don't have any expression for $f_{2}$. Two strategies are possible
here. Given the $f_{1}$ specification, one possibility is to find an
expression for $Q_{1}$, then map it through equation \eqref{gqte_eq}
to find a corresponding expression for $Q_{2}$, invert it to determine
$F_{2}$, and finally differentiate the result to get $f_{2}$. Apart
from simple situations, usually these steps (i.e. integration,
inversion and differentiation) need to be performed numerically. A
second possibility is to replace the $Y_{2}$ density in the likelihood
with its correspondent \emph{density quantile function},
$f_2(Q_2(p_{j})|\bs{\beta},\bs{\eta})$ \citep[see][]{Parzen79},
for which the next theorem provides a closed form expression. The proof
of the theorem and two additional corollaries are available in Appendix
2
, while in the next subsection we provide some further explanation on
how to compute $f_{2}$.

%
\begin{theorem}\label{theorem1}
Let $Y_1 | \bs{\eta} \sim F_1(\cdot\,;\bs{\eta})$, with $F_{1}$
having density function $f_{1}(\cdot\,;\bs{\eta})$, and assume that
\eqref{gqte_eq} holds. If, for every $0<p<1$, the vector $\bs{\beta
}$ satisfies the constraint
%
\begin{equation}\label{bsq_constraint}
X'(p,\lambda)\,\bs{\beta} \left\{\frac{d}{d\left(X(p,\lambda)\,
\bs{\beta}\right)} h^{-1}\left[X(p,\lambda)\,\bs{\beta}\right]
\right\} \leq\frac{1}{f_1\left(Q_1(p)|\bs{\eta}\right)Q_1(p)}\,,
\end{equation}
%
the density quantile function $f_2(Q_2(p)|\bs{\beta},\bs{\eta})$
for $Y_{2}$ is
\begin{equation}\label{derivsquaref3}
\textstyle f_2(Q_2(p)|\bs{\beta},\bs{\eta})=\frac{f_1\left
(Q_2(p)\,h^{-1}\left[X(p,\lambda)\,\bs{\beta}\right]|\bs{\eta
}\right) h^{-1}\left[X(p,\lambda)\,\bs{\beta}\right]}{1-f_1\left
(Q_2(p)\,h^{-1}\left[X(p,\lambda)\,\bs{\beta}\right]|\bs{\eta
}\right)\,X'(p,\lambda)\,\bs{\beta}\,Q_2(p)\left\{\frac{d}{d\left
(X(p,\lambda)\,\bs{\beta}\right)}h^{-1}\left[X(p,\lambda)\,\bs
{\beta}\right]\right\}} \,.
\end{equation}
The function $f_2(Q_2(p)|\bs{\beta},\bs{\eta})$ is a properly
defined density.
\end{theorem}

Note that $f_2$ correctly depends upon both the model parameter $\bs
{\beta}$ and the $Y_1$ parameter $\bs{\eta}$ through the $f_1$
density. The motivation for the constraint \eqref{bsq_constraint}
comes from the need to guarantee that $f_{2}$ is a non-negative
function. 
As a further remark, we observe that the term in the likelihood
involving $f_{2}(Q_{2}(p_{j}) | \bs{\beta}, \bs{\eta})$ depends
upon the observations $y_{2j}$ through the unknown quantile function
values $Q_{2}(p_{j})$.

\subsection{Details for the Computation of $f_{2}$}\label{details}
A computational drawback of our proposal is that the ``true'' values of
the percentiles $p_{j}$, i.e. those generated under the assumed model
for $Y_{2}$, should be used in the calculation of the likelihood.
Unfortunately, these are not available, because the cumulative
distribution function $F_{2}$ is not given explicitly and we cannot
find the $p_{j}$ corresponding to the observed data $y_{2j}$ as
$F_{2}(y_{2j}) = p_{j}$.

The approach we recommend to bypass this issue is to approximate the
$p_{j}$ using the procedure described in \cite{Gilchrist}, which we
summarize as follows:
\begin{enumerate}
\item Denoting with $y_{2(j)}$ the ordered observed values for $Y_{2}$,
we look for the corresponding set of ordered $p_{(j)}$ such that
$y_{2(j)} = \widehat{Q}_{2}(p_{(j)})$, where $\widehat{Q}_{2}(p)$ is
an estimate of $Q_{2}(p)$ based on the current values of the parameters
$\bs{\beta}$ and $\bs{\eta}$ (i.e. the values from the current MCMC
draw). More specifically, we find the $p_{(j)}$ using the following
procedure: suppose $p_{0}$ is the current estimate of $p$ for a given
$y$ value. Then, for a value of $p$ close to $p_{0}$, $\widehat
{Q}_{2}(p)$ can be approximated using the following Taylor series expansion
\begin{eqnarray}\label{Q2approx1}
\widehat{Q}_{2}(p) &=& \widehat{Q}_{2}(p_{0}) + \widehat
{Q}'_{2}(p_{0})(p - p_{0}) \nonumber\\
&=& \widehat{Q}_{2}(p_{0}) + \widehat{q}_{2}(p_{0})(p - p_{0}),
\nonumber
\end{eqnarray}
which, solving for $p$, gives
%
\begin{equation}\label{p2approx1}
p = p_{0} + \frac{y - \widehat{Q}_{2}(p_{0})}{\widehat{q}_{2}(p_{0})},
\end{equation}
where $\widehat{q}_{2}(p_{0})$ is the quantile density function
corresponding to $\widehat{Q}_{2}(p_{0})$ and where we used the fact
that $y = \widehat{Q}_{2}(p)$. As a starting point for $p_{(j)}$ we
use $j/(n_{2} + 1)$, $j = 1,\ldots,n_{2}$. Equation \eqref{p2approx1}
is used in an iterative fashion till the given value of $\widehat
{Q}_{2}(p)$ differs from $y$ by less than some chosen small amount (we
use $10^{-8}$).
\item Once the values of $p_{j}$ are available, we compute the
quantities $f_{2}(Q_{2}(p_{j})|\bs{\beta}, \bs{\eta})$ using
equation \eqref{derivsquaref3}. The critical issue in this step is the
calculation of the derivative $X'(p, \lambda)$. In the cases we
consider here (i.e. either a polynomial or a spline basis), the
derivative is available in closed form and so no further numerical
approximation is needed.
\end{enumerate}

Strictly speaking, the calculation of $f_{2}$ provided by the procedure
we just described is not exact but involves a numerical approximation.
We performed a detailed analysis on the goodness of this approximation
and we found that the actual and approximated $f_{2}$ values (and hence
the overall likelihood) were indistinguishable.

\subsection{Special Cases}\label{bsq_example}
We present now some special cases where an appropriate choice of the
design matrix $X(p,\lambda)$ allows to recover an exact expression for
$f_2(Q_2(p)|\bs{\beta},\bs{\eta})$ belonging to a known
distribution family. The proofs of these special cases are provided in
Appendix 3.

\noindent\textit{Case 1: $Y_1$ is Uniform and $X(p,\lambda=0)=1$.}
In this case we assume that $Y_1|\theta_1 \sim\mathcal{U}[\,0,\theta
_1]$ and choose $h(x)=x$. Then $Q_1(p)/Q_2(p) =\beta_0$ and from
\eqref{derivsquaref3_x} it follows that
\begin{equation*}
f_2(Q_2(p)|\theta_1,\beta_0) = \frac{\beta_0}{\theta_1}\,\mathbb
{I}_{\left[0,\,\theta_1/\beta_0\right]}\{Q_2(p)\}\,,
\end{equation*}
which corresponds to the density quantile function of a uniform random
variable with parameter $\theta_2=\theta_1/\beta_0$, where $\mathbb
{I}_{A}\{x\}$ denotes the indicator function taking value $1$ if $x \in
A$ and $0$ otherwise. Note that it correctly depends both upon the
parameter $\beta_0$ and the $Y_1$ parameter $\eta=\theta_1$.

\noindent\textit{Case 2: $Y_1$ is Log-normal and $X(p,\lambda=1)=[1,
\Phi^{-1}(p)]$.} We now assume $Y_1|\mu_1,\sigma_1^2 \sim\mathcal
{L}n(\mu_1,\sigma^2_1)$ and fix $h(x)=\log(x)$. It follows that
$\log\{Q_1(p)/Q_2(p)\}=\beta_0+\beta_1\,\Phi^{-1}(p)$, where $\Phi
^{-1}(p)$ is the quantile function of a standard normal random
variable. Then by \eqref{derivsquaref3_logx} we get
\begin{equation*}
f_2(Q_2(p)|\mu_1,\sigma_1^2,\beta_0,\beta_1) = \frac
{1}{Q_2(p)\sqrt{2\pi}(\sigma_1-\beta_1)}\,\exp\left\{-\frac
{\left[ \log Q_2(p) -(\mu_1-\beta_0) \right]^2}{2(\sigma_1-\beta
_1)^2}\right\}\,,
\end{equation*}
which is the density quantile function of a $\mathcal{L}n(\mu
_2,\sigma_2^2)$ random variable with $\mu_2=(\mu_1-\beta_0)$ and
$\sigma_2=(\sigma_1-\beta_1)$. Note that the density quantile
function of $Y_2$ correctly depends both upon the parameters $\bs
{\beta}=(\beta_0,\beta_1)$ and the $Y_1$ parameters $\bs{\eta
}=(\mu_1,\sigma_1^2)$. In this case the constraint \eqref
{constr_logx} simply requires that $\beta_1 \leq\sigma_1$, for every
$0 < p < 1$.

\noindent\textit{Case 3: $Y_1$ is Pareto and $X(p,\lambda=1)=[1,\log
(1-p)]$.}
Suppose $Y_1|a_1,b_1 \sim\mathcal{P}a(a_1,b_1)$ and choose $h(x)=\log
(x)$. In this case $\log\{Q_1(p)/Q_2(p)\}=\beta_0+\beta_1\log(1-p)$
and by \eqref{derivsquaref3_logx} we get
\begin{equation*}
f_2(Q_2(p)|a_1,b_1,\beta_0,\beta_1) = \frac{a_{1}}{a_{1}\beta
_{1}+1}\left(b_{1} e^{-\beta_{0}}\right)^{\frac{a_{1}}{a_{1}\beta
_{1}+1}} Q_2(p)^{-\left(\frac{a_{1}}{a_{1}\beta_{1}+1}+1\right)} \,,
\end{equation*}
which represents the density quantile function of a $\mathcal
{P}a(a_2,b_2)$ random variable with $a_2=\frac{a_{1}}{a_1\beta_1+1}$
and $b_2=b_1 e^{-\beta_0}$. The density quantile function of $Y_2$
correctly depends both upon the parameters $\bs{\beta}=(\beta
_0,\beta_1)$ and the $Y_1$ parameters $\bs{\eta}=(a_1,b_1)$. In this
case the constraint \eqref{constr_logx} requires that $\beta_1 \geq
-\frac{1}{a_1}$, for every $0 < p < 1$.

\subsection{Prior Structure and Posterior Calculation}\label{priorstructure}
The parameters $\bs{\beta}$ and $\bs{\eta}$ are assumed to be a
priori independent, that is
\[
p(\bs{\beta},\bs{\eta}|\zeta_{\bs{\beta}},\zeta_{\bs{\eta}})
= p(\bs{\beta}|\zeta_{\bs{\beta}}) \times p(\bs{\eta}|\zeta
_{\bs{\eta}})\,,
\]
where $\zeta_{\bs{\beta}}$ and $\zeta_{\bs{\eta}}$ are the prior
hyperparameters for $\bs{\beta}$ and $\bs{\eta}$ respectively. For
$p(\bs{\beta}|\zeta_{\bs{\beta}})$ we use a a multivariate normal
distribution with mean equal to the ordinary least squares (OLS)
estimate of $\bs{\beta}$ based on the model
%
\begin{equation}\label{square_ols}
h\left(\frac{y_{1(i)}}{y_{2(i)}}\right)=X(p_i,\lambda)\bs{\beta
}+\varepsilon_i\,, \quad i=1,\ldots,n,
\end{equation}
where $n=\min(n_1,n_2)$, $p_i=i/(n+1)$ and variance-covariance matrix
equal to $\sigma_{\bs{\beta}}^2 I_{(\lambda+1)}$, where $\sigma
_{\bs{\beta}}^2$ is normally fixed at a high value to induce a weakly
informative prior distribution for each $\beta_{j}$ and $I_{(\lambda
+1)}$ indicates the identity matrix with size $(\lambda+1)$. For the
prior distribution for $\bs{\eta}$, the choice clearly depends upon
the assumption made about $F_{1}$, but we suggest to use conjugate
priors. For an example see the application in Section \ref{data_analysis}.

The posterior distributions of $\bs{\beta}$ and $\bs{\eta}$ are
obtained by an MCMC simulation. In particular, we use an independent
Metropolis-Hastings algorithm with blocking over $\bs{\beta}$ and
$\bs{\eta}$ separately \citep[see][]{GilksRichSpieg,RobertCasella,OHaganForster,CarlinLouis}.
As the proposal distribution for $\bs
{\beta}$ we use a $(\lambda+1)$-dimensional $t$ distribution with
mean and scale matrix chosen to match the $\bs\beta$ OLS estimate and
variance from \eqref{square_ols}, and a small number of degrees of
freedom, usually set to $3$. As with the prior, the proposal
distribution for $\bs{\eta}$ depends upon the particular application
under investigation (see Section \ref{data_analysis} for an example).

\subsection{GQTE Estimation and Inference}\label{gqte_inference}
Once the $\bs\beta$ parameters have been estimated and the
convergence of the simulated chains has been assessed by conventional
methods \citetext{see \citealp{CarlinLouis}, or \citealp{GelmanCarlinSternRubin}},
we can obtain the posterior distribution of
the GQTE for any choice of $g(\cdot)$ as we now describe.

For each iteration $m$ of the MCMC simulation, a value $\widehat{\bs
{\beta}}^{(m)}$ for $\bs{\beta}$ is available. Using expressions
\eqref{gqte1} and \eqref{gqte2} we can obtain $\widehat
{Q}_{1}^{(m)}(p)$ and $\widehat{Q}_{2}^{(m)}(p)$ as
\[
\widehat{Q}_{1}^{(m)}(p_{2i}) = y_{2(i)} h^{-1}\left[ X\left
(p_{2i},\lambda\right)\widehat{\bs{\beta}}^{(m)} \right] \,,
\qquad i=1,\ldots,n_{2},
\]
and
\[
\widehat{Q}_{2}^{(m)}(p_{1i}) = y_{1(i)} \left\{ h^{-1}\left[ X\left
(p_{1i},\lambda\right)\widehat{\bs{\beta}}^{(m)} \right] \right\}
^{-1} \,, \qquad i=1,\ldots,n_{1},
\]
where the $y_{\ell(i)}$ are the order statistics for sample $\ell$,
while $p_{\ell i}=i/(n_{\ell}+1)$, $i=1,\ldots,n_{\ell}$, for $\ell
\in\{1,2\}$. It then follows that the $m$-th iteration value for
$\Delta_g(p)$ is given by
%
\begin{equation}\label{param_est}
\widehat{\Delta}^{(m)}_{g}(p) = g\left( \widehat{Q}_{1}^{(m)}(p)
\right)-g\left( \widehat{Q}_{2}^{(m)}(p) \right) , \qquad0<p<1,
\end{equation}
where $\widehat{Q}_{1}^{(m)}(p)$ and $\widehat{Q}_{2}^{(m)}(p)$ are
found by interpolating the estimated quantile functions $\left(
p_{2i}, \widehat{Q}_{1}^{(m)}(p_{2i}) \right)$ and $\left( p_{1i},
\widehat{Q}_{2}^{(m)}(p_{1i}) \right)$. The estimate of $\Delta
_{g}(p)$ is finally obtained through the Rao-Blackwellized estimator
%
\begin{equation}\label{param_est_fin}
\widehat{\Delta}_{g}(p) = \frac{1}{M} \sum_{m=1}^{M}\widehat
{\Delta}^{(m)}_{g}(p) \,, \qquad0<p<1,
\end{equation}
where $M$ is the total number of iterations. Since the whole posterior
distribution of $\Delta_g(p)$ is available, standard inferential
questions can be easily addressed in the usual ways.

As detailed in Section \ref{GQTE}, many interesting special cases
arise from the general definition of the GQTE. For example, if interest
lies in estimating the QTE defined in \eqref{gqte_QTE}, $g(\cdot)$
corresponds to the identity function and expression \eqref{param_est} becomes
%
\begin{equation}\label{bsq_delta_perc_formula}
\widehat{\Delta}^{(m)}(p) = \widehat{Q}_{1}^{(m)}(p) - \widehat
{Q}_{2}^{(m)}(p)\,,
\end{equation}
which can be evaluated for any value of $p \in(0,1)$.

If the focus is on the ATE, defined in \eqref{gqte_ATE}, then \eqref
{param_est} returns the estimator
\begin{multline}\label{bsq_delta_formula}
\widehat{\Delta}^{(m)} = \frac{1}{n_2}\sum^{n_2}_{i=1}{y_{2(i)}
h^{-1}\left[ X\left(p_{2i},\lambda\right)\widehat{\bs{\beta
}}^{(m)} \right]} \\
- \frac{1}{n_1}\sum^{n_1}_{i=1}{y_{1(i)} \left\{ h^{-1}\left[
X\left(p_{1i},\lambda\right)\widehat{\bs{\beta}}^{(m)} \right]
\right\}^{-1}} \,.
\end{multline}

Appendix 4
contains the details for some other cases. As a final remark, note that
an appealing feature of this approach is that the MCMC procedure needs
to be run only once to compute the difference between any measure of
the treatment effect of interest in the two groups.

\subsection{Selecting the Number of Degrees of Freedom $\lambda
$}\label{dfselection}
The choice of the number of degrees of freedom $\lambda$ to use in the
procedure above is not trivial. Many approaches can be proposed, but to
keep the computational burden manageable, we propose to elicit it by
minimizing an empirical version of the $L_1$ discrepancy measure \citep{DevLug01}
%
\begin{equation}\label{l1criterion}
D^{0}(\lambda)=\sum_{i=1}^{n} \left| f_2(Q_2(p_{i})|\bs{\beta},\bs
{\eta}) - f_2^{0}(Q_2(p_{i})) \right| ,
\end{equation}
where $f_2^0$ denotes the unknown true $Y_{2}$ density. More precisely,
we select $\lambda$ using the following procedure: for each $\lambda
\in\{ 1,\ldots,\lambda_{\max} \}$, where $\lambda_{\max}$ is the
maximum admissible value for $\lambda$, we estimate
$f_2(Q_2(p_{i})|\widehat{\bs{\beta}},\widehat{\bs{\eta}})$, where
$\widehat{\bs{\beta}}$ is equal to the OLS estimate given in \eqref
{square_ols} and $\widehat{\bs{\eta}}$ is estimated by using only
the data $(y_{11},\ldots,y_{1n_{1}})$. After replacing the true
unknown density $f_2^0$ with a kernel density estimate of the data
$(y_{21},\ldots,y_{2n_{2}})$, the value of $\lambda$ is chosen as
that minimizing the value of \eqref{l1criterion} over the set $\{
1,\ldots,\lambda_{\max} \}$. One drawback of this approach is that
it tends to select high values of $\lambda$. We provide further
discussion about this issue in Section \ref{discussion}.

\section{Simulation Study}\label{simstudy}
In this section we report the results of a simulation study we
performed to compare the finite sample properties of our method with
those of other flexible approaches.The simulation indicates that often
the GQTE procedure has lower mean squared error and a similar bias as
other flexible methods for comparing two distributions. In particular,
we contrast our proposal with the smooth quantile ratio estimation
(SQUARE) approach presented in \cite{DomCopeZeger}, and the Probit
stick-breaking process (PSBP) proposed in \cite{ChenDunson2009}. The
former is a frequentist semiparametric method while the latter is a
Bayesian nonparametric model. We now provide some details about these
two methodologies and a justification for using them.

In SQUARE it is assumed that the log quantile ratio is a smooth
function of the percentile $p$ with $\lambda$ degrees of freedom, that is
\begin{equation*}
\log\left\{ \frac{Q_{1}(p)}{Q_{2}(p)} \right\} = s(p, \gamma),
\qquad0 < p < 1.
\end{equation*}
The basic idea of smooth quantile ratio estimation is to replace the
empirical quantiles with smoother versions obtained by smoothing the
log-transformed ratio of the two quantile functions across percentiles.
SQUARE has been proposed by \cite{DomCopeZeger} as an estimator of the
mean difference between two populations with the advantage of providing
substantially lower mean squared error and bias than the sample mean
difference or the maximum likelihood estimator for log-normal
populations. To estimate $\gamma$ a $B$-fold cross-validation approach
is suggested. Finite sample inference is performed by bootstrap but
they also provide large sample results.

The PSBP is a general nonparametric Bayesian model which has been
proposed by \cite{ChenDunson2009} for estimating the conditional
distribution of a response variable given multiple predictors. More
specifically, the PSBP is a prior for an uncountable collection of
random distributions. Like for the models belonging to the class of
dependent Dirichlet processes \citep{MacEachern1999}, the PSBP main
idea is to allow for dependence across a family of related
distributions as a function of some covariates. More explicitly, the
PSBP induce dependence in the weights of the stick-breaking
representation \citep{Sethuraman1994} by replacing the
beta-distributed random variables with a probit model. This simple
change greatly enhances the flexibility of the model, thus providing an
extremely interesting extension within the framework of dependent
priors across families of probabilities measures.

We decided to compare the GQTE with these two methods because they are
both highly flexible and have been proved to perform well under a broad
set of situations. Originally the simulation also included the ANOVA
dependent Dirichlet process mixtures proposed by \cite{deIorio2004},
but we decided not to report it here because it performed poorly as
compared to the PSBP model. The reason for such inferior results
resides in the definition of the model itself, which assumes the
weights in the stick-breaking representation of the process to be
fixed, i.e. the same for the response distributions under the two
treatment conditions.

Since SQUARE produces an estimate only for the mean difference between
two populations, our simulation is restricted to this specific case. We
are aware that the results only provide a partial demonstration of the
GQTE advantages over the other methods, but we also need to stress that
the ATE is the most common measure used in practice for estimating the
extent of a treatment effect.

Our simulation framework is similar to that used in \cite{DomCopeZeger}
and includes five scenarios, which are described in
Table \ref{tabsim1} under the labels A to E. In scenarios A, B and C
the $Y_{2}$ distribution is assumed to be log-normal with parameters
$\mu_{2} = 7$ and $\sigma_{2} = 1.5$ which approximately correspond
to the sample statistics for the medical expenditures of non-diseased
subjects from the NMES data set. In scenario A, the $Y_{1}$
distribution is also log-normal but with larger values of the
parameters, namely $\mu_{1} = 7.5$ and $\sigma_{1} = 1.75$. Scenarios
B and C use a different assumption for the $Y_{1}$ distribution chosen
to represent some reasonable shapes. The next two scenarios, D and E,
compare the performances of the different methods using real data. In
particular, in scenario D the data are randomly drawn from the
distributions of nonzero Medicare expenditures for cases and controls
from the NMES data set. Finally, scenario E assumes that both
populations follow a gamma distribution with finite second moment.

%
\begin{table}[htbp]
\begin{center}
\caption{Simulation Study -- Sampling mechanisms under each simulation
scenario. In scenario D, $\widehat{F}_{g} \: (g =1, 2)$ are the
empirical cumulative distribution functions of the nonzero medical
expenditures for patients in the case and control groups from the NMES
data set, and, in scenarios B and C, $g(u)= exp\{7+1.5 \cdot\Phi
^{-1}(u)\}$. Moreover, in scenario B, $s_{B}(u) = \mathbb{I}_{(0,
1)}(u) + \mathbb{I}_{(0.9, 1)}(u)$, while in scenario C, $s_{C}(u) = 8
u (1-u)\mathbb{I}_{(0, 1)}(u)$. \newline}
\label{tabsim1}
\begin{tabular}{|c|c|c|cc|} \hline
Scenario & Population 1 & Population 2 & $n_{1}$ & $n_{2}$ \\ \hline
A & $\mathcal{L}og\mathcal{N}(7.5, 1.75)$ & $\mathcal{L}og\mathcal
{N}(7, 1.5)$ & 100 & 1000 \\
B & $u \sim\mathcal{U}nif(0, 1)$, $y_{1} = g(u) e^{s_{B}(u)}$ &
$\mathcal{L}og\mathcal{N}(7, 1.5)$ & 100 & 1000 \\
C & $u \sim\mathcal{U}nif(0, 1)$, $y_{1} = g(u) e^{s_{C}(u)}$ &
$\mathcal{L}og\mathcal{N}(7, 1.5)$ & 100 & 1000 \\
D & $\widehat{F}_{1}$ & $\widehat{F}_{2}$ & 100 & 1000 \\
E & $\mathcal{G}a(2.5, 2.5/\bar{y}_{1})$ & $\mathcal{G}a(2.5,
2.5/\bar{y}_{2})$ & 100 & 1000 \\ \hline
\end{tabular}
\end{center}
\end{table}

Under each scenario we compare the mean squared error (RMSE) and bias
(RB) in percentage relative to the sample mean difference $(\bar
{y}_{1} - \bar{y}_{2})$ for the following methods: (1) the GQTE
approach that assumes $Y_{1}$ to be log-normally distributed, (2) the
GQTE assuming $Y_{1}$ follows a gamma distribution, (3) the SQUARE
method using natural cubic splines with the number of degrees of
freedom chosen by $10-$fold cross-validation, (4) the PSBP model using
the treatment indicator as the only predictor. The GQTE estimators use
a natural cubic spline basis for the cubic-root transformed
quantile-ratio smoother with the number of degrees of freedom $\lambda
$ chosen following the procedure detailed in Section \ref
{dfselection}. The RMSE is computed by $[ \{ \mbox{\normalfont
{MSE}}(\bar{y}_{1}-\bar{y}_{2})-\mbox{\normalfont{MSE}}(\widehat
{\Delta}) \} / \mbox{\normalfont{MSE}}(\bar{y}_{1}-\bar{y}_{2})]
\times100$, while the RB is defined as $[ \{ \mathbb{E}(\widehat
{\Delta})-\Delta\} / \Delta] \times100$. Note that positive values
for the RMSE imply a better performance for the estimator as compared
to the sample mean difference.

The results for $100$ generated data sets for each scenario are
reported in Table \ref{tabsim2}. We considered only the case of
unbalanced samples with $n_{1} = 100$ and $n_{2} = 1000$ because
typically it represents a more critical situation to deal with in
practice. These results show that the GQTE has a smaller mean squared
error in most of the scenarios considered.

%
\begin{table}[htbp]
\begin{center}
\caption{Simulation Study -- Results from $100$ replicate datasets.
RMSE is the mean squared error relative to $(\bar{y}_{1}-\bar
{y}_{2})$ in percentage defined by $[ \{ \mbox{\normalfont{MSE}}(\bar
{y}_{1}-\bar{y}_{2})-\mbox{\normalfont{MSE}}(\widehat{\Delta}) \}
/ \mbox{\normalfont{MSE}}(\bar{y}_{1}-\bar{y}_{2})] \times100$,
and RB is the bias relative to $(\bar{y}_{1}-\bar{y}_{2})$ in
percentage defined by $[ \{ \mathbb{E}(\widehat{\Delta})-\Delta\} /
\Delta] \times100$, under the data generation mechanisms described in
Table \ref{tabsim1}. The splines degrees of freedom $\lambda$ for the
GQTE approach are chosen using the heuristic algorithm described in
Section \ref{dfselection} while for SQUARE we use 10-fold
cross-validation.}\label{tabsim2}\vspace*{9pt}
\footnotesize
\begin{tabular}{|@{\ }c|cc|cc|cc|cc|cc|} \hline
& \multicolumn{2}{c |}{Scenario A} & \multicolumn{2}{c |}{Scenario B}
& \multicolumn{2}{c |}{Scenario C} & \multicolumn{2}{c |}{Scenario D}
& \multicolumn{2}{c |}{Scenario E} \\
& RMSE & RB & RMSE & RB & RMSE & RB & RMSE & RB & RMSE & RB \\ \hline
GQTE ($\mathcal{L}og\mathcal{N}$) & 39 & -29 & -26 & -23 & 3 & 13 &
-26 & 18 & -1 & -2 \\
GQTE ($\mathcal{G}amma$) & 33 & -16 & 1 & -9 & 48 & 7 & 25 & 0 & 0 &
-2 \\
SQUARE & 35 & -6 & -7 & -5 & 1 & 13 & 25 & 3 & 0 & -1 \\
PSBP & 1 & -6 & 6 & -10 & -8 & 9 & -13 & 6 & -3 & -3 \\ \hline
MSE$(\bar{y}_1-\bar{y}_2)$ & \multicolumn{2}{c |}{2952} &
\multicolumn{2}{c |}{5992} & \multicolumn{2}{c |}{1051} &
\multicolumn{2}{c |}{2100} & \multicolumn{2}{c |}{753} \\
$\Delta$ & \multicolumn{2}{c |}{4982} & \multicolumn{2}{c |}{15225}
& \multicolumn{2}{c |}{5244} & \multicolumn{2}{c |}{7144} &
\multicolumn{2}{c |}{7144} \\ \hline
\end{tabular}
\end{center}
\end{table}

In scenario A, where both the populations are log-normal, the GQTE
assuming $Y_{1}$ is log-normally distributed performs around $40\%$
better than $(\bar{y}_{1} - \bar{y}_{2})$, slightly better than
SQUARE, even if somewhat biased. This result is superior to that of
PSBP, which performs approximately as well as the sample mean
difference. 
In scenario B, the PSBP provides the best result with a mean square
error which is $6\%$ smaller than $(\bar{y}_{1} - \bar{y}_{2})$,
followed by the GQTE with gamma distributed $Y_{1}$. In scenarios C, D
and E the GQTE with gamma distributed $Y_{1}$ outperforms both the PSBP
and SQUARE. More specifically, in scenario C the GQTE provides a mean
square error that is approximately $50\%$ smaller than $(\bar{y}_{1} -
\bar{y}_{2})$. This is also the least biased result. In scenarios D
and E, the GQTE approach provides comparable results as those provided
by SQUARE.

\section{Application: Medical Costs for Smoking Attributable
Diseases}\label{data_analysis}
As an illustration, we apply the GQTE approach to the NMES data, where
the distributions of $Y_{1}$ (the cases) and $Y_{2}$ (the controls) are
highly right-skewed. For this reason, we decide to use $h(x)=\log(x)$.
We show that having a smoking attributable disease induces both a
location and scale shift in the medical expenditure distribution as
compared to that for non-affected subjects, but with a thinning of the
corresponding distribution's tails.

\subsection{Data Description}
The data used in the following analysis is taken from the National
Medical Expenditure Survey (NMES) and have been previously studied by
other authors \citep[for example][]{DomCopeZeger}. It provides data on
annual medical expenditures, disease status, age, race, socio-economic
factors, and critical information on health risk behaviors such as
smoking, for a representative sample of U.S. non-institutionalized
adults \citep{NMES}. NMES data derive from the 1987 wave. 
In the data set used here a total of 9,416 individuals are available.
Table \ref{NMES} briefly summarizes the data set (numbers in
parentheses represent the percentage of subjects with non-zero expenditures).

%
\begin{table}[htbp]
\caption{Disease cases and controls for smokers (current or former)
and for non-smokers. Numbers within parentheses represent the
percentage of people in that cell with non-zero expenditures.}\label{NMES}
\begin{center}
\begin{tabular}{|l|c|c|c|} \hline
& smokers & non smokers & Total \\ \hline
cases & 165 (62\%) & 23 (70\%) & 188 (63\%) \\ \hline
controls & 4,682 (21\%) & 4,546 (28\%) & 9,228 (25\%) \\ \hline
Total & 4,847 (22\%) & 4,569 (28\%) & 9,416 (25\%) \\ \hline
\end{tabular}
\end{center}
\end{table}

We consider as cases ($Y_{1}$) those individuals who are affected by
smoking diseases, namely lung cancer and chronic obstructive pulmonary
disease, while the controls ($Y_{2}$) are persons without a major
smoking attributable disease.

In the following analyses we consider only the non-zero costs paid for
each hospitalization by diseased and non-diseased subjects.

\begin{figure}
\includegraphics{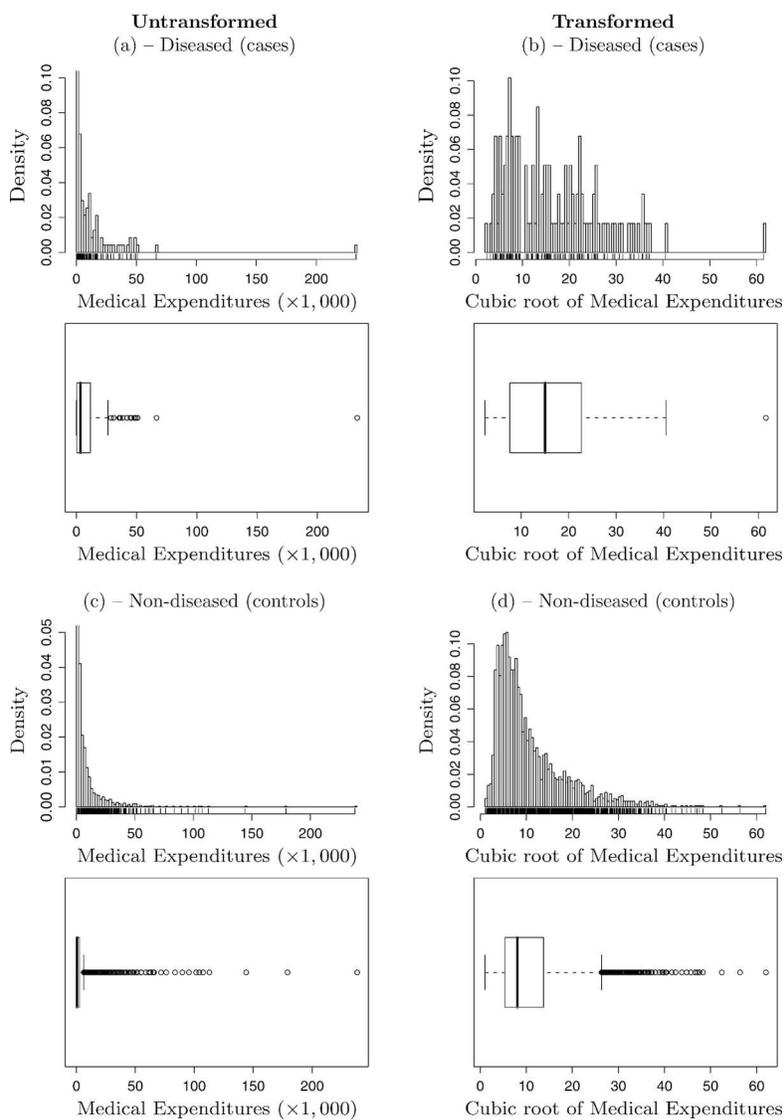}
\caption{Histograms and boxplots of positive medical
expenditures for hospitalizations regarding smoking attributable
diseases (lung cancer and coronary obstructive pulmonary disease) from
the 1987 National Medicare Expenditure Survey (for clarity of
exposition, the histogram of the original expenditures has been
truncated at the top).}\label{data_summary}
\end{figure}

Figures \ref{data_summary}(a) and (c) show the histograms and boxplots
for the medical costs of the cases and controls. Both the distributions
are highly right-skewed, with the cases sample which is much smaller
than the controls one ($118$ vs. $2,262$). Table \ref{table_summary}
contains some high-order sample quantiles for the two groups which
confirm the heavier tails of the cases costs distribution. However,
note also that the controls sample has a higher maximum cost.

%

%
\begin{table}[htbp]
\caption{Summary of the NMES data set: high-order quantiles of
non-zero medical expenditures for cases and controls.}
\begin{center}
{\scriptsize
\begin{tabular}{|c|c|c|c|c|c|c|c|c|}
\hline
Quantile order & 75 & 90 & 95 & 99 & 99.9 & 100 \\
\hline
Quantile for cases (\$) & 11,525.17 & 29,439.96 & 49,595.77 & 63,886.05
& 213,567.69 & 233,047.63 \\
\hline
Quantile for controls (\$) & 2,600.00 & 9,799.664 & 30,625.206 &
49,771.60 & 135,896.07 & 238,185.94 \\
\hline
\end{tabular}
}
\end{center}
\label{table_summary}
\end{table}

Figures \ref{data_summary}(b) and (d) show the histograms and boxplots
for the cubic root transformed data. The need for such a transformation
derives from the particular choice we make regarding the cases
distribution (see next subsection) and is not a general requirement of
our approach. Moreover, the use of this transformation does not alter
in any way the results and the conclusions we draw, hence in the
following we systematically refer to the transformed data. At any rate,
note also that, even after the transformation, the outcome
distributions still present heavy right tails. This conclusion
motivates the use of $h(x)=\log(x)$.

\begin{figure}[h]
\includegraphics{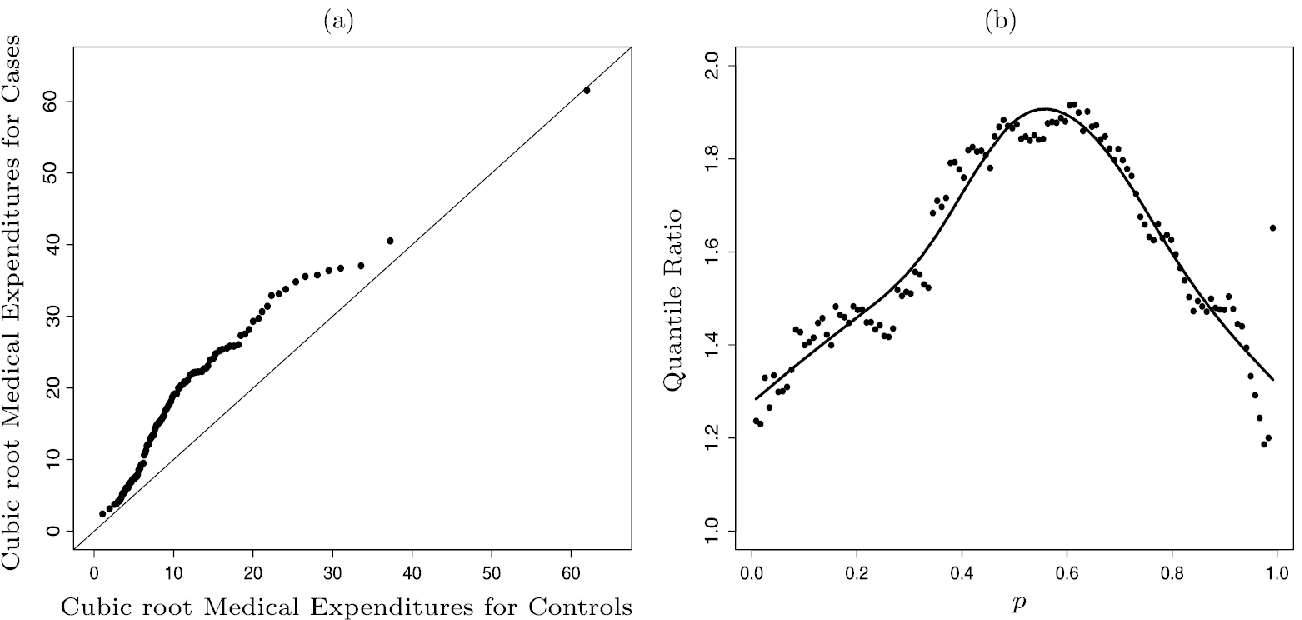}
\caption{(a) Q-Q plot of cubic root transformed non-zero medical
expenditures. (b)~Quantile ratio across percentiles with a fitted
natural cubic spline.}\label{square_qqplot_cubic}
\end{figure}

In Figure \ref{square_qqplot_cubic}(a) we report the Q-Q plot for the
(cubic root transformed) NMES data. We can identify a non-linear smooth
relationship between the cases and controls medical expenditures. Panel
(b) of the same picture, which shows the quantile ratio as a function
of the percentile $p$, confirms these findings.

%

\subsection{Model Assumptions and Tuning Parameters}
In this application we assume that $Y_1|\bs{\pi},\theta\sim\mathcal
{GSM}\left(\bs{\pi},\theta| J\right)$, a particular mixture of
gamma distributions with density
\begin{equation*}\label{gsm}
f(y|\pi_1,\ldots,\pi_J,\theta)=\sum_{j=1}^{J}{\pi_j \frac{\theta
^{j}}{\Gamma(j)}y^{j-1}e^{-\theta y}}\,,
\end{equation*}
where the mixing occurs over the shape parameters and where $J$, the
number of components, is fixed a priori. First introduced in \citet{VentGSM},
it has been explicitly developed as a model for right-skewed
distributions and its parameterization allows to create a convenient
and flexible method characterized by a single scale parameter for all
the gamma components, plus the ordinary set of mixture weights. We use
conjugate priors $\theta\sim\mathcal{G}a(\alpha,\delta)$ and $\bs
{\pi} \sim\mathcal{D}_{J}\left(\frac{1}{J},\ldots,\frac
{1}{J}\right)$ for the shared scale parameter and the mixture weights
respectively. The number of mixture components is fixed at $J=40$,
while the $\theta$ hyperparameters are set to $\alpha=845$ and
$\delta=1,300$ \citep[for more information on the elicitation of
these priors see][Section~2.3]{VentGSM}.

The initial values of the Metropolis-Hastings algorithm are chosen as
follows: the $\bs{\beta}$ chain is started from its OLS estimate, as
discussed in \eqref{square_ols}, while for $\bs{\eta}=(\theta,\bs
{\pi})$ we first get a preliminary estimate (with $5,000$ iterations)
using the approach described in \citet{VentGSM}, and then we fix their
starting values to the corresponding estimated posterior averages.

We run the MCMC algorithm for $1,000,000$ iterations plus $200,000$
iterations as burn-in. Such a large number of iterations is necessary
because the model, being quite complicated, has shown a slow
convergence behavior of the chains.

\subsection{Results}
The selection procedure described in Subsection \ref{dfselection} for
the number of degrees of freedom suggested a value of $\lambda$ equal
to $6$, which can be considered fairly satisfactory from a visual
inspection of the scatterplot (see Figure \ref{square_qqplot_cubic}(b)).

The acceptance rates for the MCMC posterior simulations are relatively
small, being around $0.5\%$ for $\bs{\beta}$, $25\%$ for $\theta$
and $1.6\%$ for $\bs{\pi}$. Despite that, we do not consider these
results as problematic since the chain is moving in a high-dimensional
space ($(\lambda+ 1) + (J + 1) = 48$ dimensions), which necessarily
slows down the convergence process. This is the main reason why we
decide to run the simulation for a longer time. However, the results of
the analysis presented below indicate that convergence was attained. We
made other attempts with simpler (but less flexible) specifications of
the $Y_1$ distribution, which showed a more conventional behavior of
the acceptance rates.\looseness=1

\begin{figure}[h]
\includegraphics{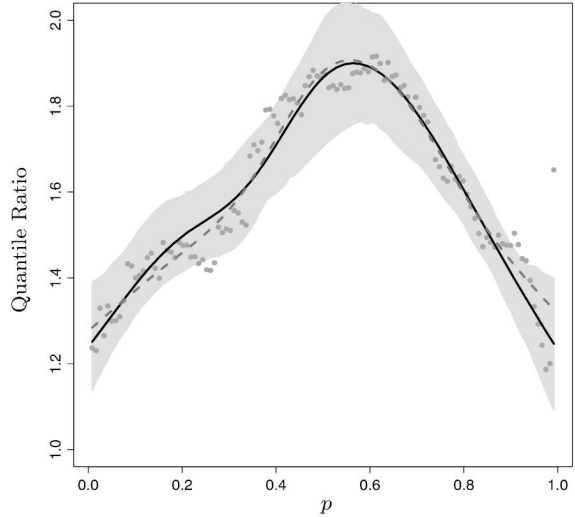}
\caption{Fitted values of the estimated model \eqref{gqte_eq}. The
solid line represents the estimated pointwise posterior means, while
the shaded area corresponds to their pointwise $95\%$ credible
intervals. The dashed line corresponds to the OLS fit for the same
data, as described in \eqref{square_ols}.}\label{bsq_fit}
\end{figure}

Figure \ref{bsq_fit} shows the fitted values for the estimated model
\eqref{gqte_eq}. The gray dots represent the quantile ratio for the
transformed data as a function of the percentile $p$. The solid line
illustrates the estimated posterior mean of the quantile ratio, while
the dashed one represents its OLS estimate (the same line as in Figure
\ref{square_qqplot_cubic}(b)). The shaded area gives the credible
bands for the estimated posterior means, showing a fairly low amount of
uncertainty around the estimates. Moreover, from the picture we can see
that our model is less sensitive to extreme observations, especially in
the right tails of the distributions.\looseness=1

%

\begin{figure}[h]
\includegraphics{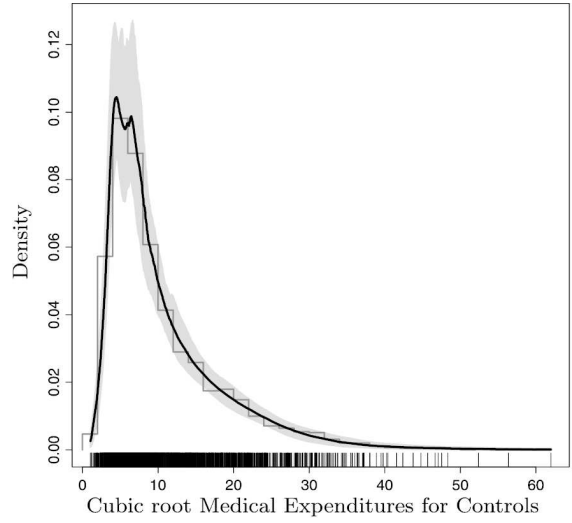}
\caption{Estimated $Y_2$ density. The solid line represents the
estimated pointwise posterior means, while the shaded area shows the
corresponding $95\%$ credible intervals. The thinner dark gray line
depicts the data histogram.}\label{est_f2}
\end{figure}

In Figure \ref{est_f2} we report the estimated $Y_2$ density. It is
possible to ascertain a quite good fit. In the display, together with
the $f_2(Q_2(p_{2i})|\theta,\bs{\pi},\bs{\beta})$ posterior mean,
we put the corresponding $95\%$ credible bands and the histogram of the data.

%

%
\begin{figure}[t!]
\includegraphics{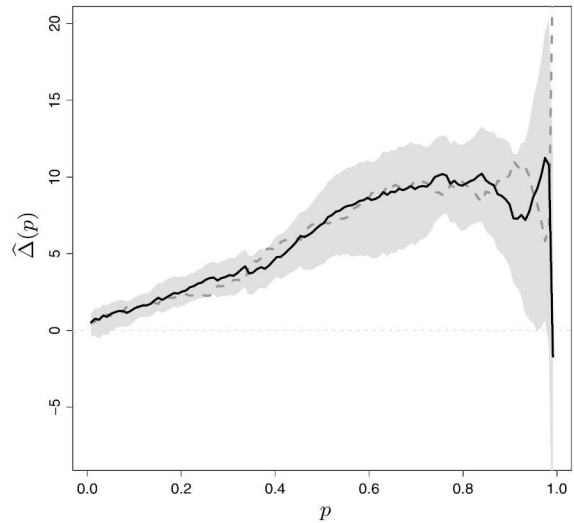}
\caption{Estimated Quantile Treatment Effect (QTE), defined as $\Delta
(p)=Q_1(p)-Q_2(p)$. The solid line reports the estimated pointwise
posterior means, the shaded area gives the corresponding $95\%$
credible intervals, while the dashed line shows the sample quantile
differences. Data are cubic-root transformed.}\label{delta_id}
\end{figure}

We now describe the results for the GQTE $\Delta_g(p)$ introduced in
Section \ref{estimation} for some choices of the function $g(\cdot)$.
We start from the QTE, denoted as $\widehat{\Delta}(p)$ in \eqref
{bsq_delta_perc_formula}, whose estimate is shown in Figure \ref
{delta_id}. The solid line represents the posterior mean of the medical
costs QTE between cases and controls and the gray area is the
corresponding $95\%$ credible interval, while the dashed line portrays
the sample quantile differences. We can see that the distribution of
the medical expenditures for subjects with smoking attributable
diseases is always above that of those without smoking-related
diseases. However, a much larger variability results in estimating the
difference for the very extreme quantiles. This behavior is not too
surprising since the two samples become very sparse as the medical
expenditures become bigger (see the boxplots in Figure~\ref{data_summary}).

\begin{figure}[t]
\includegraphics{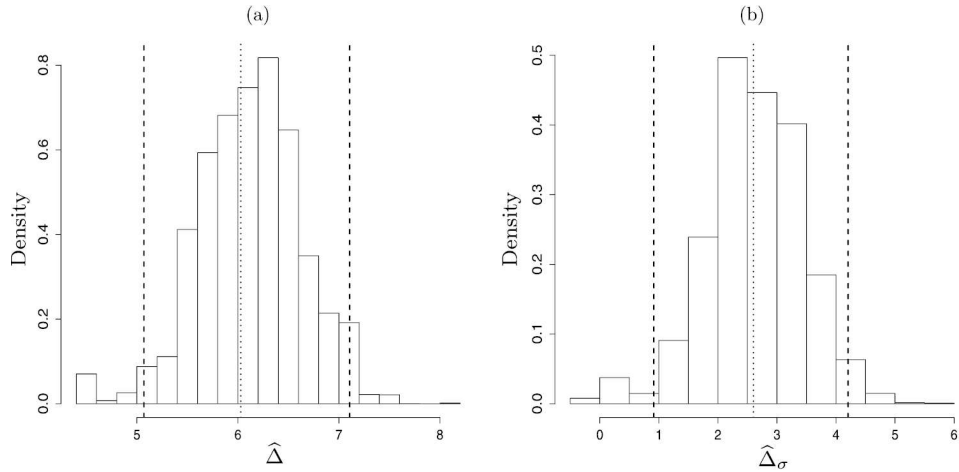}
\caption{(a) Estimated posterior distribution of the Average Treatment
Effect (ATE) between cases and controls medical expenditures (cubic
root transformed), as defined in \eqref{gqte_ATE}. The vertical dashed
lines represent the $95\%$ credible interval, while the dotted line is
the sample mean difference. (b) Estimated posterior distribution of the
standard deviation difference between cases and controls medical
expenditures (cubic root transformed), as defined in Section \ref
{GQTE}. The vertical dashed lines represent the $95\%$ credible
interval, while the dotted line depicts the sample standard deviation
difference.}\label{delta_mu}
\end{figure}

Figures \ref{delta_mu}(a) and (b) contain the posterior distributions
of the ATE and the standard deviation difference, as defined in Section
\ref{GQTE}. The sample mean and standard deviation differences are
depicted in the two plots with a vertical dotted line, while the $95\%$
credible intervals are indicated using dashed lines. The ATE estimated
posterior mean (on the log scale) is equal to $6.1127$, while the
estimated posterior mean of the standard deviation difference is
$2.6275$. These results prove that having a smoking attributable
disease has a significant negative impact on both the location and
scale of the single hospitalization medical cost distribution. After
re-transforming the estimated quantiles on the original scale%
, we get an estimated posterior mean for the ATE between diseased and
non-diseased subjects equal to $\$6,244.10$.

%

\begin{figure}[h]
\includegraphics{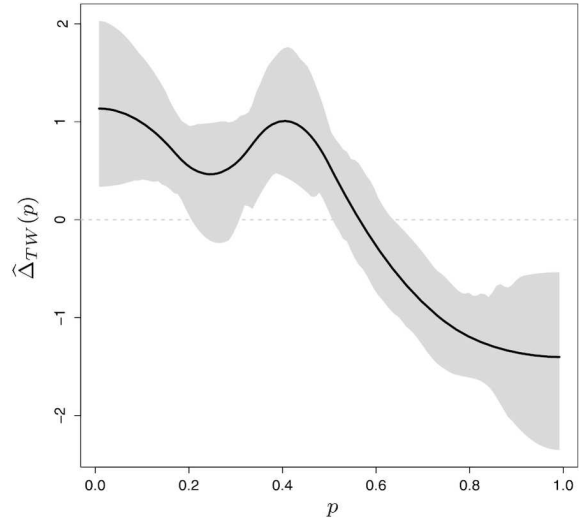}
\caption{Estimated tailweight difference $\Delta
_{TW}(p)=TW_1(p)-TW_2(p)$ between cases and controls, as defined in
\eqref{gqte_tw}. The solid line represents the estimated pointwise
posterior means, while the shaded area shows the corresponding $95\%$
credible intervals.}\label{delta_tw}
\end{figure}

Finally, Figure \ref{delta_tw} shows the impact of the treatment
variable (i.e., having or not a smoking attributable disease) on the
tailweight functions $TW(p)$, defined in \eqref{gqte_tw}, of the two
populations. The tails of the medical costs distribution for the
diseased subjects tend to be heavier than those of the non-diseased
ones for values of $p$ up to approximately $0.6$, but the situation is
inverted as we move to consider higher percentiles. Hence, while the
fact of being affected by smoking attributable diseases tends to
increase both the average and the variance of the medical expenditures
distribution, we have found that the opposite occurs to the tail
probabilities, that is, to the chances of incurring very high medical
costs in a single hospitalization.

As a last comment, we would like to remark on the explicit choice we
made to exclude the observations with null medical costs. We took this
decision because the inclusion of this further feature of the data
requires the extension of our approach to a two-part modeling framework
\citep{Mullahy,CameronTrivedi}, which doesn't appear to be
straightforward in our context.

%


\section{Discussion}\label{discussion}
In this paper we have introduced a new parameter, the GQTE, for
assessing the effect of a binary covariate on a response and a novel
methodology to estimate it. The GQTE generalizes the most common
approaches available in the literature, that is, the well-known average
treatment effect (ATE) and the quantile treatment effect (QTE), since
it allows to evaluate the effect of a treatment on any arbitrary
characteristic of the outcome's distributions under the two treatment
conditions.

To estimate the GQTE we have proposed a Bayesian procedure, where we
assume that a monotone transformation of the quantile ratio is modeled
as a smooth function of the percentiles. This assumption allows to
increase efficiency by borrowing information across the two groups. The
idea of quantile ratio smoothing has first been introduced by \citet{DomCopeZeger}.
In the present work we extended that proposal in
several ways: 1) we let the link between the quantile ratio and the
percentiles be general and application-specific, allowing to take into
account the tail heaviness of the distributions involved in the
analysis; 2) we derive a closed form expression for the model
likelihood; 3) our methodology is not limited to the mean difference
between the treated and the controls, but provides a comprehensive
assessment of the treatment effect; 4) finally, we embed the whole
estimation process within a Bayesian framework allowing to make
inference on the GQTE $\Delta_{g}(p)$ for any choice of the function
$g(\cdot)$, and for both symmetric and highly skewed outcomes.

The GQTE is a marginal measure in the sense that it provides an
estimate of the treatment effect over an entire population. In the
econometrics literature this kind of approach is usually termed the
\textit{unconditional} QTE \citep{Firpo2007,FrolichMelly2008} in
contrast with the \textit{conditional} QTE, where the treatment effect
is determined separately for different combinations of a set of
covariates \citep{KoenkerBassett,Koenker,AngristPischke}. The
inclusion of covariates can improve the efficiency of an estimator even
when the primary goal of the analysis is a marginal effect.
Accordingly, methods have been proposed to extract marginal quantiles
from estimates of conditional quantiles \citep{MachadoMata,FrolichMelly2008}.
A challenge in extending our approach along these
lines is the lack of an ``iterated expectation'' result\footnote{While
for a standard linear model, in fact, the assumption $E(Y_{i}|X_{i}) =
X_{i}'\bs{\beta}$ does imply $E(Y_{i}) = E(X_{i})'\bs{\beta}$, the
same conclusion doesn't hold for the conditional quantiles.} for the
quantiles \citep[see for example][Chapter 7]{AngristPischke}.

To further clarify our goals, we want to stress that in this paper no
particular emphasis has been placed on the causality issues that
naturally comes into play when the objective is the estimation of a
treatment effect \citep[see for example][]{Rosenbaum2002,Rosenbaum2010,Rubin2006,AngristPischke}.
More precisely, our intent
here is solely to provide a general measure of the effect of a binary
treatment on a response variable, together with a flexible approach to
estimate it.

We compared the performance of our estimation approach with other
highly flexible methods in a simulation study for the mean difference
between two populations. Our study revealed that the GQTE performs
generally better than the other competing estimators at least in
estimating the mean difference.


We have applied our methodology to the NMES data set to assess the
effect of being affected by smoking attributable diseases on the single
hospitalization medical costs distribution. We have found that having
these diseases increases the average medical bill amount as well as its
variability in the population, while it reduces the probability of
incurring higher bills.

Our approach can be extended in various directions. The most promising
research question we can see involves taking into account individual
level characteristics in measuring the effect of a treatment. In our
context, this would involve the estimation of a conditional version of
$\Delta_{g}(p)$, something like $\Delta_{g}(p | \bs{x})=g\left(
Q_{1}(p | \bs{x}) \right) - g\left( Q_{2}(p | \bs{x}) \right)$.
The clear advantage of including covariates would be an increase in the
efficiency of the estimates \citep{FrolichMelly2008}. To control for
systematic differences in covariates between two populations, a common
strategy is to group units into subclasses based on covariate values,
for example using propensity score matching, and then to apply our
method within strata of propensity scores \citep{Rosenbaum2002,Rosenbaum2010},
as implemented for example in \citet{DomZeger}.

Currently we are considering only a binary treatment effect, so another
important line of research is the extension of the methods to
categorical ordinal and to continuous treatments.

A further direction for future research concerns the choice of the
number of degrees of freedom $\lambda$. In this paper we adopted the
simple approach of choosing $\lambda$ by minimizing an empirical
version of the $L_{1}$ distance between the $Y_{2}$ density estimate
and its kernel density estimate (see Subsection \ref{dfselection}).
More structured solutions can obviously be\vadjust{\eject} considered. A natural
extension would allow $\lambda$ to be a random quantity to be
estimated together with all the other parameters using a
trans-dimensional MCMC approach, like for example the reversible jump
algorithm \citep{green1995}. While this solution would allow to take
into account also the uncertainty connected to the a priori ignorance
about the $\lambda$ value, the consequence would be a dramatic
increase in the computational workload of the estimation algorithm.

\section*{Acknowledgements}
The research of Dominici was supported by Award Number R01ES012054
(Statistical Methods for Population Health Research on Chemical
Mixtures) from NIH/NIEHS, Award Numbers R83622 (Statistical Models for
Estimating the Health Impact of Air Quality Regulations) and RD83241701
(Estimation of the Risks to Human Health of PM and PM Components) from
EPA, Award Number 4909-RFA11-1/12-3 (Causal Inference Methods for
Estimating Long Term Health Effects of Air Quality Regulations) from
HEI and Award Number K18 HS021991 (A Translational Framework for
Methodological Rigor to Improve Patient Centered Outcomes in End of
Life Cancer Research) from AHRQ. The content is solely the
responsibility of the authors and does not necessarily represent the
official views of the above Institutions.

\section*{Appendix 1: Additional GQTE Examples}\label{appendixA}
Together with the cases presented in Section \ref{GQTE}, many other
less conventional measures of the difference between two distributions
can be obtained by properly choosing the $g(\cdot)$ function in the
GQTE definition. For example, by choosing $g(x)=\int{x^r}\,dp$ we
obtain the difference between the population $r$-th moments
%
\begin{equation}\label{gqte_rmoment}
\Delta_{\mu^r}=\int_0^1{Q_1(p)^r}\,dp-\int_0^1{Q_2(p)^r}\,dp.
\end{equation}
Using the fact that for a random variable $Y$ with expected value $\mu
$, variance $\sigma^2$ and quantile function $Q(p)$ it holds that (see
\citealp{Gilchrist} or \citealp{Shorack})
\begin{equation*}
\sigma^2=\int_0^1{\left[ Q(p) - \mu\right]^2}\,dp=\int_0^1{
Q(p)^{2} }\,dp - \mu^{2} ,
\end{equation*}
by suitably choosing the $g(\cdot)$ function, we recover the
difference between the two population variances as
%
\begin{eqnarray}\label{gqte_variance}
\Delta_{\sigma^2} &=& \left[\int_0^1{Q_1(p)^2}\,dp-\left(\int
_0^1{Q_1(p)}\,dp\right)^2\right]-\left[\int_0^1{Q_2(p)^2}\,dp-\left
(\int_0^1{Q_2(p)}\,dp\right)^2\right] \nonumber\\
&=& \left[\int_0^1{Q_1(p)^2}\,dp-\int_0^1{Q_2(p)^2}\,dp\right
]-\left[\left(\int_0^1{Q_1(p)}\,dp\right)^2-\left(\int
_0^1{Q_2(p)}\,dp\right)^2\right] \nonumber\\
&=& \Delta_{\mu^2}-\left(\mu_1^2-\mu_2^2\right),
\end{eqnarray}

However, the cases encompassed by the GQTE include many other
quantile-based indexes that are less frequently used in the literature,
like the inter-$p$-range $ipr(p)=Q(1-p)-Q(p)$, or the skewness-ratio
$sr(p)=[Q(1-p)-Q(0.5)]/[Q(0.5)-Q(p)]$, $0<p<1$, which provide robust
measures of the scale and shape of a distribution \citep[for a list of
these indexes see][]{Gilchrist,Shorack,Parzen04,WangSerfling,BrysHubertStruyf}.
A quantity of particular interest to economists is
the difference between inter-decile ratios, defined as
\[
\frac{Q_{1}(0.9)}{Q_{1}(0.1)}-\frac{Q_{2}(0.9)}{Q_{2}(0.1)},
\]
which is commonly used to measure the inequality in a population \citep
[see][]{FrolichMelly2008}. The previous quantity can be easily
generalized as follows
%
\begin{equation}\label{inter_p_range}
\Delta_{IR}(p) = \frac{Q_{1}(1-p)}{Q_{1}(p)}-\frac{Q_{2}(1-p)}{Q_{2}(p)},
\end{equation}
for any $0 < p < 0.5$. Notice that all these indexes are obtainable
from the general definition \eqref{GQTE_def} by properly choosing the
function $g(\cdot)$.

As a last example, we consider a further GQTE special case that is
based on the so called \emph{tailweight function} defined as
\begin{equation*}\label{tailweight}
TW(p) = \frac{q(p)}{Q(p)} \equiv\frac{d}{dp}\log Q(p) \,, \qquad
0<p<1\,,
\end{equation*}
which is used to quantify the probability allocated in the tails of a
distribution. One can compute the difference between the tailweight
functions for two populations by choosing the logarithmic derivative of
the quantile function as the $g(\cdot)$ functional in \eqref
{GQTE_def}, that is
%
\begin{eqnarray}\label{gqte_tw}
\Delta_{TW}(p) &=& TW_{1}(p)-TW_{2}(p) \nonumber\\
&=& \frac{d}{dp}\log Q_{1}(p) - \frac{d}{dp}\log Q_{2}(p) \nonumber\\
&=& \frac{d}{dp} \left[ \log\left( \frac{Q_{1}(p)}{Q_{2}(p)}
\right) \right].
\end{eqnarray}
If $\Delta_{TW}(p) \geq0$, we can conclude that the treatment is
causing a thickening of the $Y_{1}$ distribution tails as compared to
those of $Y_{2}$ if $\Delta_{TW}(p) \geq0$. 
Finally, note that, thanks to the equivariance property of the
quantiles, \eqref{gqte_tw} can be written also as
%
\begin{eqnarray}\label{gqte_tw_2}
\Delta_{TW}(p) &=& \frac{d}{dp}\log Q_{1}(p) - \frac{d}{dp}\log
Q_{2}(p) \nonumber\\
&=& \frac{d}{dp}Q_{1,\log}(p) - \frac{d}{dp} Q_{2,\log}(p)
\nonumber\\
&=& \frac{d}{dp}\left[ Q_{1,\log}(p) - Q_{2,\log}(p) \right]
\nonumber\\
&=& \frac{d}{dp} \Delta_{\log}(p) \,,
\end{eqnarray}
where $Q_{\ell,\log}$, $\ell=1,2$, indicates the quantile of the
log-transformed data and $\Delta_{\log}(p)$ denotes the parameter
\eqref{gqte_QTE} calculated on the quantiles of the log-transformed data.

\section*{Appendix 2: Proof of Theorem \ref{theorem1} and
Corollaries}\label{appendixB}
\begin{proof}[Proof of Theorem 1]\label{proof_th1}
Differentiate \eqref{gqte1} with respect to $p$ to get
%
\begin{equation}\label{derivsquare}
q_1(p)=q_2(p)\,h^{-1}\left[X(p,\lambda)\,\bs{\beta}\right
]+X'(p,\lambda)\,\bs{\beta}\,Q_2(p)\left\{\frac{d}{d\left
(X(p,\lambda)\,\bs{\beta}\right)}h^{-1}\left[X(p,\lambda)\,\bs
{\beta}\right]\right\},
\end{equation}
where $q_{\ell}(p)=dQ_{\ell}(p)/dp$ denotes the so called \textit
{quantile density function} for the population $\ell=\{1,2\}$, while
$X'(p,\lambda)$ corresponds to the derivative of $X(p,\lambda)$,
$0<p<1$ (properly resized because a constant is normally included in
the design matrix $X(p,\lambda)$).

Apply now to both $q_1(p)$ and $q_2(p)$ the following relationship
between the quantile density and density quantile functions \citep[see
for example][]{Gilchrist,Parzen79,Parzen04}
%
\begin{equation}\label{quandens}
f(Q(p))\,q(p)=1\,,
\end{equation}
to get the expression
\begin{multline}\label{derivsquaref}
\frac{1}{f_1(Q_1(p)|\bs{\eta})}=\frac{1}{f_2(Q_2(p)|\bs{\beta
},\bs{\eta})}\,\,h^{-1}\left[X(p,\lambda)\,\bs{\beta}\right] \\
+X'(p,\lambda)\,\bs{\beta}\,Q_2(p) \left\{\frac{d}{d\left
(X(p,\lambda)\,\bs{\beta}\right)}h^{-1}\left[X(p,\lambda)\,\bs
{\beta}\right]\right\},
\end{multline}
and hence
%
\begin{equation}\label{derivsquaref2}
\textstyle f_2(Q_2(p)|\bs{\beta},\bs{\eta})=\frac{f_1\left
(Q_1(p)|\bs{\eta}\right) h^{-1}\left[X(p,\lambda)\,\bs{\beta
}\right]}{1-f_1\left(Q_1(p)|\bs{\eta}\right)\,X'(p,\lambda)\,\bs
{\beta}\,Q_2(p)\left\{\frac{d}{d\left(X(p,\lambda)\,\bs{\beta
}\right)}h^{-1}\left[X(p,\lambda)\,\bs{\beta}\right]\right\}} \,.
\end{equation}
Finally, substituting \eqref{gqte1} in place of $Q_{1}(p)$ proves the
main statement.

Moreover, $f_2(Q_2(p)|\bs{\beta},\bs{\eta})$ is a proper density
function because:
\begin{itemize}
\item$f_2(Q_2(p)|\bs{\beta},\bs{\eta}) \geq0$, for any $0<p<1$;
since $f_1\left(Q_1(p)|\bs{\eta}\right) \geq0$ and $h^{-1}\left
[X(p,\lambda)\,\bs{\beta}\right] > 0$ because, as assumed in \eqref
{gqte_eq}, it is the ratio of two positive quantile functions, this
fact can be proved by showing that the denominator of \eqref
{derivsquaref2} is nonnegative which is ensured by the constraint
\eqref{bsq_constraint}.
\item$\int_{0}^{1}f_2(Q_2(p)|\bs{\beta},\bs{\eta}) q_2(p) dp =
1$, which is true because $f_2(Q_2(p)|\bs{\beta},\bs{\eta}) q_2(p)
= 1$ by construction.\qedhere
\end{itemize}
\end{proof}

A couple of immediate consequences of Theorem \ref{theorem1} regard
two cases that occur frequently in practice. We provide the details
about these situations in the next two corollaries.

%
\begin{coroll}
Let the same assumptions of Theorem \ref{theorem1} hold. Suppose
additionally that $h(x)=x$. If for every $0<p<1$ the vector $\bs{\beta
}$ satisfies the constraint
%
\begin{equation}\label{constr_x}
\frac{X'(p,\lambda)\,\bs{\beta}}{X(p,\lambda)\,\bs{\beta}} \leq
\frac{1}{f_1\left(Q_1(p)\right)Q_1(p)}\,,
\end{equation}
then the density quantile function $f_2(Q_2(p)|\bs{\beta},\bs{\eta
})$ for $Y_{2}$ is
%
\begin{equation}\label{derivsquaref3_x}
f_2(Q_2(p)|\bs{\beta},\bs{\eta})=\frac{f_1\left(Q_2(p)\,
X(p,\lambda)\,\bs{\beta}|\bs{\eta}\right)X(p,\lambda)\,\bs
{\beta}}{1 -f_1\left(Q_2(p)\,X(p,\lambda)\,\bs{\beta}|\bs{\eta
}\right)X'(p,\lambda)\,\bs{\beta}\,Q_2(p)}\,.
\end{equation}
\end{coroll}

%
\begin{coroll}
Let the same assumptions of Theorem \ref{theorem1} hold. Suppose
additionally that $h(x)=\log(x)$. If for every $0<p<1$ the vector $\bs
{\beta}$ satisfies the constraint
%
\begin{equation}\label{constr_logx}
X'(p,\lambda)\,\bs{\beta} \leq\frac{1}{f_1\left(Q_1(p)\right
)Q_1(p)}\,,
\end{equation}
then the density quantile function $f_2(Q_2(p)|\bs{\beta},\bs{\eta
})$ for $Y_{2}$ is given by
%
\begin{equation}\label{derivsquaref3_logx}
f_2(Q_2(p)|\bs{\beta},\bs{\eta})=\frac{f_1\left(Q_2(p)\,e^{\,
X(p,\lambda)\,\bs{\beta}}|\bs{\eta}\right)}{\,e^{-\,X(p,\lambda
)\,\bs{\beta}} -f_1\left(Q_2(p)\,e^{\,X(p,\lambda)\,\bs{\beta
}}|\bs{\eta}\right)X'(p,\lambda)\,\bs{\beta}\,Q_2(p)}\,.
\end{equation}
\end{coroll}

Note that in these two situations, the general constraint \eqref
{bsq_constraint} reduces to a linear constraint on $\bs\beta$.

\section*{Appendix 3: Proofs of the Special Cases}\label{appendixC}
\noindent\textit{Case 1: $Y_1$ is Uniform and $X(p,\lambda=0)=1$.}
Here $Y_1|\theta_1 \sim\mathcal{U}[\,0,\theta_1]$ and $h(x)=x$. In
this case $Q_1(p)/Q_2(p) =\beta_0$. Hence the density, distribution
and quantile functions of $Y_1$ are respectively
\begin{eqnarray*}
f_1(y_1|\theta_1) &=& \frac{1}{\theta_1}\,\mathbb{I}_{[0,\,\theta
_1]}\{y_1\} \\
F_1(y_1|\theta_1) &=& \frac{y_1}{\theta_1} \\
Q_1(p|\theta_1) &=& \theta_1 p\,, \qquad\qquad0<p<1\,.
\end{eqnarray*}
From \eqref{derivsquaref3_x} it follows that
\begin{eqnarray*}
f_2(Q_2(p)|\theta_1,\beta_0) &=& \frac{1}{\theta_1}\,\mathbb
{I}_{[0,\,\theta_1]}\{Q_2(p)\,\beta_0\} \, \beta_0 \\
&=& \frac{\beta_0}{\theta_1}\,\mathbb{I}_{\left[0,\,\theta
_1/\beta_0\right]}\{Q_2(p)\}\,,
\end{eqnarray*}
which is the density quantile function of a $\mathcal{U}[0,\theta_2]$
random variable with $\theta_2=\theta_1/\beta_0$.

\noindent\textit{Case 2: $Y_1$ is Log-normal and $X(p,\lambda=1)=[1,
\Phi^{-1}(p)]$.} Assume $Y_1|\mu_1,\sigma_1^2 \sim\mathcal{L}n(\mu
_1,\sigma^2_1)$ and $h(x)=\log(x)$. In this case $\log\{
Q_1(p)/Q_2(p)\}=\beta_0+\beta_1\,\Phi^{-1}(p)$, where $\Phi
^{-1}(p)$ is the quantile function of a standard normal random
variable. The density, distribution and quantile functions of $Y_1$ are
given by
\begin{eqnarray*}
f_1(y_1|\mu_1,\sigma_1^2) &=& \frac{1}{y_1\sqrt{2\pi}\sigma
_1}\exp\left\{-\frac{\left( \log y_1-\mu_1 \right)^2}{2\sigma
_1^2}\right\} \\
F_1(y_1|\mu_1,\sigma_1^2) &=& \Phi\left(\frac{\log y_1 - \mu
_1}{\sigma_1}\right) \\
Q_1(p|\mu_1,\sigma_1^2) &=& \exp\left\{\,\mu_1+\sigma_1\Phi
^{-1}(p)\right\}\,, \qquad\quad0<p<1\,.
\end{eqnarray*}
Then by \eqref{derivsquaref3_logx} it follows
{\small
\begin{eqnarray*}
f_2(Q_2(p)|\mu_1,\sigma_1^2,\beta_0,\beta_1) &=& \frac{\frac
{1}{\sqrt{2\pi}\sigma_1}\,\frac{\exp\left\{-\frac{\left( \mu
_1+\sigma_1\Phi^{-1}(p)-\mu_1 \right)^2}{2\sigma_1^2}\right\}
}{\exp\left\{\,\mu_1+\sigma_1\Phi^{-1}(p)\right\}}}{\exp\left\{
-\beta_0-\beta_1\Phi^{-1}(p)\right\}\left[ 1-\frac{1}{\sqrt{2\pi
}\sigma_1}\,\frac{\exp\left\{-\frac{\left( \mu_1+\sigma_1\Phi
^{-1}(p)-\mu_1 \right)^2}{2\sigma_1^2}\right\}}{\exp\left\{\,\mu
_1+\sigma_1\Phi^{-1}(p)\right\}} \right.} \\
&& \qquad\qquad\frac{\phantom{1}}{\left. \times\beta_1 \frac
{1}{\frac{1}{\sqrt{2\pi}}\exp\left\{-\frac{\left[\Phi
^{-1}(p)\right]^2}{2}\right\}} \, \exp\left\{\,\mu_1+\sigma_1\Phi
^{-1}(p)\right\} \right]} \\
&=& \frac{1}{\sqrt{2\pi}(\sigma_1-\beta_1)}\frac{\exp\left\{
-\frac{\left[ \Phi^{-1}(p) \right]^2}{2}\right\}}{\exp\left\{\,
(\mu_1-\beta_0)+(\sigma_1-\beta_1)\Phi^{-1}(p)\right\}} \\
&=& \frac{1}{\sqrt{2\pi}(\sigma_1-\beta_1)}\,\frac{\exp\left\{
-\frac{\left[ (\mu_1-\beta_0)+(\sigma_1-\beta_1)\Phi
^{-1}(p)-(\mu_1-\beta_0) \right]^2}{2(\sigma_1-\beta_1)^2}\right\}
}{\exp\left\{\,(\mu_1-\beta_0)+(\sigma_1-\beta_1)\Phi
^{-1}(p)\right\}} \\
&=& \frac{1}{Q_2(p)\sqrt{2\pi}(\sigma_1-\beta_1)}\,\exp\left\{
-\frac{\left[ \log Q_2(p) -(\mu_1-\beta_0) \right]^2}{2(\sigma
_1-\beta_1)^2}\right\}\,,
\end{eqnarray*}
}%
which is the density quantile function of a $\mathcal{L}n(\mu
_2,\sigma_2^2)$ random variable with $\mu_2=(\mu_1-\beta_0)$ and
$\sigma_2=(\sigma_1-\beta_1)$.

\noindent\textit{Case 3: $Y_1$ is Pareto and $X(p,\lambda=1)=[1,\log
(1-p)]$.}
Now $Y_1|a_1,b_1 \sim\mathcal{P}a(a_1,b_1)$ and $h(x)=\log(x)$. In
this case $\log\{Q_1(p)/Q_2(p)\}=\beta_0+\beta_1\log(1-p)$. The
density, distribution and quantile functions of $Y_1$ are given by
\begin{eqnarray*}
f_1(y_1|a_1,b_1) &=& a_{1} b_{1}^{a_{1}} y_{1}^{-(a_{1}+1)} \\
F_1(y_1|a_1,b_1) &=& 1-b_{1}^{a_{1}} y_{1}^{-a_{1}} \\
Q_1(p | a_1, b_1) &=& b_{1} (1-p)^{-\frac{1}{a_1}} \,, \qquad0<p<1\,.
\end{eqnarray*}
Then \eqref{derivsquaref3_logx} implies
{\small
\begin{eqnarray*}
f_2(Q_2(p)|a_1,b_1,\beta_0,\beta_1) &=& \textstyle\frac
{a_{1}b_{1}^{a_{1}}\left[b_{1}(1-p)^{-\frac{1}{a_{1}}}\right
]^{-(a_{1}+1)}}{\exp\left\{-\beta_{0}-\beta_{1}\log(1-p)\right\}
\left\{1+a_{1}b_{1}^{a_{1}}\left[b_{1}(1-p)^{-\frac{1}{a_{1}}}\right
]^{-(a_{1}+1)}\frac{\beta_{1}}{1-p}b_{1}(1-p)^{-\frac
{1}{a_{1}}}\right\}} \\
&=& \frac{a_{1}}{a_{1}\beta_{1}+1}\left(b_{1} e^{-\beta_{0}}\right
)^{-1} (1-p)^
{\frac{a_{1}\beta_{1}+1}{a_{1}}+1} \\
&=& \frac{a_{1}}{a_{1}\beta_{1}+1}\left(b_{1} e^{-\beta_{0}}\right
)^{\frac{a_{1}}{a_{1}\beta_{1}+1}} Q_2(p)^{-\left(\frac
{a_{1}}{a_{1}\beta_{1}+1}+1\right)} \,,
\end{eqnarray*}
}%
which is the density quantile function of a $\mathcal{P}a(a_2,b_2)$
random variable with $a_2=\frac{a_{1}}{a_1\beta_1+1}$ and $b_2=b_1
e^{-\beta_0}$.

\section*{Appendix 4: Details About the Estimation of Other
Cases}\label{appendixD}
One can estimate the impact of a binary treatment on the $r$-th
moments, denoted as $\Delta_{\mu^r}$ in \eqref{gqte_rmoment}, by computing
\begin{multline}\label{bsq_delta_r_formula}
\widehat{\Delta}^{(m)}_{\mu^r} = \frac{1}{n_2}\sum
^{n_2}_{i=1}{\left\{ y_{2(i)} h^{-1}\left[ X\left(p_{2i},\lambda
\right)\widehat{\bs{\beta}}^{(m)} \right] \right\}^r} \\
- \frac{1}{n_1}\sum^{n_1}_{i=1}{\left[ y_{1(i)} \left\{ h^{-1}\left
[ X\left(p_{1i},\lambda\right)\widehat{\bs{\beta}}^{(m)} \right]
\right\}^{-1} \right]^r} \,.
\end{multline}
The last expression allows to estimate the treatment effect on the
population variances, defined in \eqref{gqte_variance}, which is given by
\begin{multline}\label{bsq_delta_var_formula}
\widehat{\Delta}^{(m)}_{\sigma^2} = \widehat{\Delta}^{(m)}_{\mu
^2} - \left\{ \left( \frac{1}{n_2}\sum^{n_2}_{i=1}{y_{2(i)}
h^{-1}\left[ X\left(p_{2i},\lambda\right)\widehat{\bs{\beta
}}^{(m)} \right]} \right)^2 \right. \\
\left. - \left( \frac{1}{n_1}\sum^{n_1}_{i=1}{y_{1(i)} \left\{
h^{-1}\left[ X\left(p_{1i},\lambda\right)\widehat{\bs{\beta
}}^{(m)} \right] \right\}^{-1}} \right)^2 \right\} \,.
\end{multline}
%

As a concluding example, the effect of a binary treatment on the
tailweight functions of two distributions, introduced in \eqref
{gqte_tw}, can be obtained by first computing the posterior draws
%
\begin{equation}\label{bsq_delta_tw_formula_est}
\widehat{\Delta}^{(m)}_{TW}(p) = \frac{d}{dp} \left\{ \log\left(
h^{-1}\left[ X(p,\lambda) \widehat{\bs{\beta}}^{(m)} \right]
\right) \right\} ,
\end{equation}
and then by applying \eqref{param_est_fin}. When $h(x)=\log(x)$,
\eqref{bsq_delta_tw_formula_est} becomes
%
\begin{equation}
\widehat{\Delta}^{(m)}_{TW}(p) = X'(p,\lambda) \widehat{\bs{\beta}}^{(m)},
\end{equation}
and the estimate of $\Delta_{TW}(p)$ is
%
\begin{eqnarray}
\widehat{\Delta}_{TW}(p) &=& \frac{1}{M}\sum_{m=1}^{M}{X'(p,\lambda
) \widehat{\bs{\beta}}^{(m)}} \nonumber\\
&=& X'(p,\lambda) \left( \frac{1}{M}\sum_{m=1}^{M}{\widehat{\bs
{\beta}}^{(m)}} \right) \nonumber\\
&=& X'(p,\lambda) \widehat{\bs{\beta}},
\end{eqnarray}
with $\widehat{\bs{\beta}}=\frac{1}{M}\sum_{m=1}^{M}{\widehat{\bs
{\beta}}^{(m)}}$, the posterior mean estimate of $\bs{\beta}$.



%


\end{document}